\documentclass[10pt]{amsart}

\usepackage{amsmath}
\usepackage{amsthm}
\usepackage{amsfonts}
\usepackage{amssymb,latexsym}
\usepackage[all]{xy}
\usepackage{graphicx}
\usepackage[mathscr]{eucal}
\usepackage{verbatim}
\usepackage{hyperref}

\theoremstyle{plain}
    \newtheorem{thm}{Theorem}[section]
    \newtheorem{prop}[thm]{Proposition}

    \newtheorem{cor}[thm]{Corollary}

\theoremstyle{definition}

\theoremstyle{remark}

\numberwithin{equation}{section}

\newcommand{\rar}{\ensuremath{\rightarrow}}
\newcommand{\lrar}{\ensuremath{\longrightarrow}}

\newcommand{\la}{\langle}
\newcommand{\ra}{\rangle}

\newcommand{\G}{\mathcal G}

\newcommand{\stk}[1]{\stackrel{#1}{\rightarrow}}
\newcommand{\lstk}[1]{\stackrel{#1}{\longrightarrow}}

\newcommand{\Gal}{\textup{Gal}}

\newcommand{\ftwo}{\mathbb{F}_2}
\newcommand{\fp}{\mathbb{F}_p}
\newcommand{\sep}{\text{sep}}
\newcommand{\ann}{\text{ann}}
\renewcommand{\cor}{\text{cor}}

\def\res{\operatorname{res}\nolimits}

\begin{document}

% Article information
\title[Absolute Galois groups and the Bloch-Kato conjecture]{Absolute Galois groups viewed from small quotients
and the Bloch-Kato conjecture}
\date{\today}

% author one information
\author{Sunil K. Chebolu}
\address{Department of Mathematics \\
Illinois State University \\
Normal, IL 61790, USA} \email{schebol@ilstu.edu}
%\thanks{}

% author two information
\author{J\'{a}n Min\'{a}\v{c} $^{*}$}
\address{Department of Mathematics\\
University of Western Ontario\\
London, ON N6A 5B7, Canada}
\email{minac@uwo.ca}
\thanks{J\'{a}n Min\'{a}\v{c} was supported in part by Natural Sciences and Engineering Research Council of Canada
grant R3276A01}

% AMS information

% Abstract
\begin{abstract}
In this paper we concentrate on the relations between the structure of small Galois groups,
arithmetic of fields, Bloch-Kato conjecture, and Galois groups of maximal pro-$p$-quotients
of absolute Galois groups.
\end{abstract}

\dedicatory{Dedicated to Professor Paulo Ribenboim with admiration, respect
and friendship on the occasion of his 80th birthday.}

\maketitle
\thispagestyle{empty}

\tableofcontents

% Body of the paper
\section{Introduction}
The second author fondly remembers how a number of years ago Paulo Ribenboim
helped him to escape to the West and immediately upon his arrival welcomed him
with beautiful lectures on the Galois group of the Pythagorean closure of
$\mathbb{Q}$. Paulo Ribenboim's lectures, writings, and research have
influenced us strongly, and in particular this paper reflects his influence on the choice of topics
and our way of thinking about them. The paper is a selective survey of results on small quotients of absolute Galois groups and their relations with the Bloch-Kato conjecture. It is by no
means a comprehensive historical survey. Instead, it focuses only on some
selective topics from the work of the authors and their collaborators.

The main idea of our paper, and a key point we want to illustrate, is that
already  relatively small quotients of absolute Galois groups encode
substantial information about them. Absolute Galois groups of fields play a
central role in arithmetic, geometry, and topology. Yet these profinite
groups are mysterious and not much is known about the  fundamental problem
of characterizing  absolute Galois groups among all profinite groups. Therefore,
it is natural to investigate the maximal pro-$p$-quotients of absolute
Galois groups for a fixed prime $p$ which are in general much simpler than the absolute Galois
groups.  But even these pro-$p$-quotients are
quite mysterious and one would like to find more manageable, yet interesting,
quotients of maximal pro-$p$-quotients of absolute Galois groups. One
extremely interesting family of such quotients are the $W$-groups defined below in Section 6.
These are pro-$p$-groups of nilpotent class at most $2$ and they have exponent
dividing $p^2$. So they are rather simple groups in comparison with the
maximal pro-$p$-quotients of absolute Galois groups. Nevertheless, when the
primitive $p$-th root of unity is contained in the base field, these groups
carry complete information about the entire Galois cohomology with
$\fp$-coefficients of the absolute Galois groups. This is a consequence of
the Bloch-Kato conjecture which was proved recently by M. Rost and V.
Voevodsky with C. Weibel's patch. (See \cite{Vo1} and  \cite{Vo2} for an overview of the
proof, and other references \cite{HW, MVW, Ro, Ro2, Su, Vo3, W1, W2, W3} for
the foundation and completion of the proof,
and some further exposition.) Therefore, it is clear that the $W$-groups
 of fields are good candidates for thorough investigation.

The plan of this survey paper is as follows. We begin with \v{S}afarevi\v{c}'s
early work as a motivation for studying the Galois module structure of $p$-power
classes of fields. This is interesting even in the case when we consider just
cyclic extensions of degree $p$. We show that the answer in this case already
leads to a description of relatively small pro-$p$-quotients of absolute
Galois groups called $T$-groups. We then describe some recent results on the
Galois module structure of Galois cohomology. After that we give a description
of $W$-groups, their relationship with Witt rings of quadratic forms, Galois
cohomology, valuations, and the structure of maximal pro-$p$-extensions. Here the
Bloch-Kato conjecture plays a very important role, especially in the case
$n=2$ which was established by Merkurjev and Suslin almost 30 years ago. At
the conclusion of the paper we touch upon our recent work with D. Benson and J. Swallow
in progress whose goal is to provide a refinement of the Bloch-Kato conjecture with
group cohomology, combinatorics, and Galois theoretic consequences.

\section{Some work of \v{S}afarevi\v{c}}

I.R. \v{S}afarevi\v{c} initiated the very interesting program of studying
Galois groups of maximal $p$-extensions. Let $F$ be a field and $F_{sep}$ the separable closure
of $F$. Given a prime number $p$,  let $F(p)$ be the maximal $p$-extension over $F$. Thus $F(p)$ is the
compositum in $F_{sep}$  of all Galois
extensions $K/F$ which have degree a power of $p$. Let $G_F(p) :=
\Gal(F(p)/F)$  be the Galois group of the maximal $p$-extension.
 In 1947 \v{S}afarevi\v{c} showed that if $\mathbb{Q}_p \subseteq F$, $[F :
\mathbb{Q}_p] < \infty$, and if $F$ does not contain a primitive $p$th root of
unity, then $G_F(p)$ is a free pro-$p$-group on $[F: \mathbb{Q}_p] + 1$
generators. The key part of \v{S}afarevi\v{c}'s argument was the determination of the
number $d(H)$ of minimal generators for all open subgroups $H$ of $G_F(p)$.
This number $d(H)$ is equal to $\dim_{\fp} H/\Phi(H)$, where $\Phi(H) = H^p[H,
H]$ is the Frattini subgroup of $H$. By local class field theory it follows
that $d(H) = \dim_{\fp} K^*/{K^*}^p$ where $K$ is the fixed field of $H$ and $K^{*} = K  \backslash \{0\}$
is the multiplicative group of $K$. We
can calculate $d(H)$ explicitly in terms of invariants of the extension $K/F$. It turns out that
\[\dim_{\fp} K^*/{K^*}^p - 1 = [K:F] (\dim_{\fp} F^*/{F^*}^p - 1 ).\]
Thus
\[ d(H) - 1 = (d(G_F(p))- 1)[K:F] \]
 and therefore the number of
generators $d(H)$, for each $H$ open subgroup of $G_F(p)$, grow in the same way as if
$G_F(p)$ were a free pro-$p$-group. It was the insight of \v{S}afarevi\v{c}
that in fact this property is enough to prove that $G_F(p)$ is a free
pro-$p$-group. \v{S}afarevi\v{c} did not use the convenient language of
profinite groups as this terminology was not available at that time.
Similarly, the language of Galois cohomology appeared much later, and all
became wide-spread only after the appearance of Serre's lecture notes  in
1964; see \cite{serre-GaloisCohomology} for the latest edition of Serre's book.
 In particular, now we can rewrite the above equation as
\[ \dim_{\fp} H^1(H, \fp) - 1 = (\dim_{\fp} H^1(G_F(p), \fp) -1) [G_F(p): H]. \]

Let $G$ be a pro-$p$-group with $\dim_{\fp} H^r(G, \fp) < \infty$ for $1 \le r
\le n$. Following H. Koch  we can set the $n$th partial Euler-Poincar\'{e}
characteristic $\chi_n(G)$ as
\[ \chi_n(G) = \sum_{r = 0}^n (-1)^r \dim_{\fp} H^r(G, \fp). \]
Koch proved that  if $W$ is a system of open subgroups of $G$ which form a
neighbourhood basis at 1 and $\chi_n(U) = [G: U] \chi_n(G)$ for all $U$ in
$W$, then the cohomological dimension $cd(G)$ of $G$ is at most $n$. Because
$cd(G) = 1$ if and only if $G$ is a free pro-$p$-group,  we see that Koch's
criterion for $cd(G) \le n$ generalises \v{S}afarevi\v{c}'s criterion for
$G_F(p)$ being a free pro-$p$-group.

\section{The Bloch-Kato conjecture}
We briefly recall the Bloch-Kato conjecture for the novice. This conjecture will be used in this article at various places.
Let $F$ be a field that has a primitive $p$-th root $\zeta_p$ of unity. Consider the Kummer sequence
\[ 1 \lrar \mu_{p} \lrar F_{\sep}^{*} \lstk{x \rar x^{p}} F_{{\sep}} \lrar 1\]
of modules over the absolute Galois group $G_{F}$, where $\mu_{p}$ denotes the group of
$p$-th roots of unity.
The boundary map
\[ H^{0}(G_{F}, F_{\sep}^{*}) = F^{*} \rar H^{1}(G_{F}, \fp)\]
 in the induced long exact sequence in Galois cohomology extends naturally to a map
\[ T(F^{*}) \rar H^{*}(G_{F}, \fp),\]
where $T(F^{*})$ is the tensor algebra on $F^{*}$ and $H^{*}(G_{F}, \fp)$ is the Galois
cohomology ring of $G_{F}$.  Bass and Tate verified that the Steinberg relations $(a) \cup (1-a) = 0$
for $a \ne 0, 1$ hold in $H^{*}(G_{F}, \fp)$.  Milnor $K$-theory $K_{*}(F)$ is a graded ring obtained by taking
 the quotient of $T(F^{*})$ by the
graded two-sided ideal generated by the elements $a \otimes (1-a)$,  $a \in F^{*} \backslash \{1\}$. Thus we get a map, known as the norm-residue map,
\[ \eta \colon K_*(F)/pK_*(F) \lrar H^*(F, \mathbb{F}_p )\]
from the reduced Milnor $K$-theory to the Galois cohomology of $F$.
The Bloch-Kato conjecture claims that the norm-residue  map $\eta$
is an isomorphism. The case $p = 2$ was implicitly conjectured
by J. Milnor in 1970; see \cite{Mil-1970}. The Milnor conjecture was eventually proved by Voevodsky \cite{Vo1, Vo2}. For this spectacular work Voevodsky was awarded a Fields medal in 2002. His work used some
sophisticated machinery such as motivic cohomology operations and the development  of
$\mathbb{A}^1$ stable homotopy theory.
The proof of the Bloch-Kato conjecture is even more subtle. Although Voevodsky announced a proof of the Bloch-Kato conjecture
in 2003, not until September 2007 were all of the details for the  Rost-Voevodsky proof
made available by Voevodsky, Rost and Weibel; see \cite{Ro, Ro2, Su, Vo3, W1, W2, W3}. The work on the proof of the Bloch-Kato conjecture and the resulting theorem has already had tremendous impact on  contemporary
mathematics and is expected to have an even  broader impact in the coming years. Note that the Bloch-Kato conjecture gives
a presentation of the rather mysterious Galois cohomology  $H^*(F, \mathbb{F}_p)$ by generators  and relations.
In particular, it tells us that $H^*(F, \mathbb{F}_p)$ is generated by one dimensional classes.

\section{Classical Hilbert 90 and absolute Galois groups}

\v{S}afarevi\v{c}'s approach to $G_F(p)$ made clear that the $p$-th power class
group $F^*/{F^*}^p$ is a very useful and fundamental object to study. In 1960
Faddeev began to study the Galois module structure of $p$-th power classes of
cyclic extensions of local fields, and during the mid 1960s he and Borevi\v{c}
established the Galois module structure of $p$-th power class groups of local
fields using basic arithmetic invariants attached to these extensions. (See
\cite{Faddeev, Borevic}).

In the theory of quadratic forms, the exact sequence of square class
groups associated with the quadratic extension $K = F(\sqrt{a})$, $a \in
F^*/{F^*}^2$,  $\text{char}(F) \ne 2$ has been playing an important role.
This sequence is
\[ 1 \rar \{{F^*}^2, a {F^*}^2 \} \rar F^*/{F^*}^2 \stk{i_{F/K}} K^*/{K^*}^2 \stk{N_{K/F}} F^*/{F^*}^2
\stk{\epsilon} B(F),\]
where $i_{F/K}$ is the map induced by inclusion map $F \rar K$ and $N_{K/F}
\colon K^*/{K^*}^2 \rar F^*/{F^*}^2$ is the map induced by the norm map $K^*
\rar F^*$, and $\epsilon$ is the homomorphism from $F^*/{F^*}^2$ to the Brauer
group $B(F)$ defined by $b{F^*}^2 \rar [\left( \frac{a, b}{F}\right)]$ in
$B(F)$, where $[(\frac{a, b}{F})]$ is the class of Quaternion algebra
\[(\frac{a, b}{F}) = \{ f_0 + f_1 i + f_2 j + f_3 ij  | f_l \in F, i^2 = j^2 = -1, ij = -ji \} \]
in the Brauer group of $F$. (See \cite[Pg 200, Thm 3.2]{Lam}.) Observe that this
sequence completely determines the size of $K^*/{K^*}^2$ provided we know the
size $N_{K/F}(K^*)/{F^*}^2$. In fact, this sequence determines the structure of
$K^*/{K^*}^2$ as an $\ftwo[\Gal(K/F)]$-module provided we know
$N_{K/F}(K^*)/{F^*}^2$. Therefore it is desirable to extend the  work of Borevi\v{c}
and Fadeev from the case of local fields to general fields. Borevi\v{c}, Fadeev
and \v{S}afarevi\v{c} used local class field theory to establish their results
in the case of local fields. However, in \cite{minac-swallow} it was observed that  it is
possible to determine the structure of the $\fp[\text{Gal}(K/F)]$-module $K^*/{K^*}^p$ in
the case when a primitive $p$-th root of unity is contained in $F$ using just
Hilbert 90 in place of local class field theory.  In \cite{MSS}, the work of \cite{minac-swallow}
was extended to all cyclic extensions $K/F$ of degree $p^n$ with no
restriction on the base field. The remarkable feature of the final result is
that $K^*/{K^*}^p$  can be written as a sum of cyclic modules over
$\fp[\Gal(K/F)]$ of dimensions over $\fp$ all powers of $p$ with the possible
single cyclic module exception of dimension $p^m + 1$, $0 \le m \le n-1$.

As an example we formulate here the main result of \cite{MSS} when the exceptional
summand does not occur. For the more complicated case when the exceptional
summand does occur, we refer the reader to \cite{MSS}. Let $F$ be any field.
Consider a cyclic Galois extension $K/F$ with Galois group $ G = \Gal(K/F) =
\la \sigma \ra$, a cyclic group of order $p^n$. We set $J = K^*/{K^*}^p$. Then
$J$ has the obvious $\fp G$-module structure.

\begin{thm}
Assume that $F$ does not contain any primitive $p$-th root of unity. Then the $\fp
G$-module $J$ decomposes as
\[ J \cong Y_0 \oplus Y_1 \oplus \cdots \oplus Y_n,\]
where $Y_i$ is a direct sum of cyclic $\fp G$-modules of dimension $p^i$.
\end{thm}

Moreover, the multiplicity of cyclic summands of dimension $p^i$ in $Y_i$ is
completely determined by the filtration of $F^*$ by norm groups
\[N_{K/F}(K^*) \subset N_{K_{n-1}/F}(K_{n-1}^*) \subset \cdots \subset F^* \]
where, for each $i = 0, 1, \dots n$, $K_i$ is the unique subfield of $K$ which
has dimension $p^i$ over $F$. Indeed if
\[Y_i = \oplus_{I} C_{p^i} \]
where $C_{p^i}$ is the cyclic $\fp G$-module of dimension $p^i$, then the
cardinality of $I$ is just
\[ \text{dim}_{\fp} N_{K_i/F}(K_i^*)/N_{K_{i+1}/F}(K_{i+1}^*).\]
Further set $[K_i^*] = K_i K^*/ {K^*}^p$ and set $H_i = \text{Gal}(K/K_{i})$.
Then we have Galois descent in the sense that $[K_i^*] = J^{H_i}$ -- the fixed
elements of $J$ under the action of $H_i$.

The results in the case when $F$ contains a primitive $p$-th root of unity are
similar, but technically more challenging due to the occurence of the
exceptional module mentioned earlier.

These results are remarkable because of the absence of summands of dimensions
not equal to a power of $p$. (Except in the case of exceptional summands which
have dimension $p^m + 1$.)  These results severely restrict possible small
quotients of absolute Galois groups. We shall illustrate this by describing a
result from a recent paper \cite{BLMS}.

We call  a pro-$p$-group $R$ elementary abelian if it has the form $R =
\prod_{I} C_{p}$, where $C_{p}$ is a cyclic group of order $p$, and $I$ is
some possibly infinite index set. We say that  a pro-$p$-group $G$ is a
$T$-group if $G$ contains a maximal closed subgroup $N$, $N \ne G$, of
exponent dividing $p$.  Then $N$ is a normal subgroup of $G$ and the factor
group $G/N$ acts naturally on $N$ via conjugation.  Furthermore, the subgroup
$N$ is uniquely determined by $G$ provided that $G$ is neither an elementary
abelian group of order greater than $p$ nor  the direct product of an
elementary abelian group and a non-abelian group of order $p^{3}$ of exponent
$p$ if $p > 2$, and the dihedral group of order $8$ if $p =2$. Given any
profinite group $A$ with a closed normal subgroup $B$ of index $p$, the factor
group $A/B^{p}[B,B]$ is a $T$-group. Now suppose that $E/F$ is a cyclic field
extension of degree $p$. We define the $T$-group of $E/F$ to be $T_{{E/F}}:=
G_{F}/G_{E}^{p}[G_{E}, G_{E}]$, where $G_{F}$ and $G_{E}$ are absolute Galois
groups of $F$ and $E$ respectively. (For the benefit of  topologists, in order
to avoid a possible confusion, we remark that the name ``T-group" is not
motivated by Kazhdan's property (T), and in fact we do not know of any
connection between T-groups and groups with property (T) in Kazhdan's sense.)
We shall now classify those $T$-groups which are realisable as $T_{{E/F}}$ for
fields containing a primitive $p$-th root of unity.

In order to illustrate the restrictions on those $T$-groups which are realisable as $T_{{E/F}}$ for fields which
contain a primitive $p$-th root of unity we shall introduce a simple set of invariants which determine $T$-groups
up to isomorphism. We shall then see that for $p > 2$ we obtain restrictions on possible invariants of $T$-groups
which are $T_{{E/F}}$ for suitable $E/F$. On the other hand there are no restrictions in the case $p=2$.  The
proof of these statements can be found in \cite{BLMS}.

For a pro-$p$-group $A$, denote $Z(A)$ its center and $Z(A)[p]$ the elements of $Z(A)$ of order dividing $p$.
Let $A_{(n)}$ be the $n$-th group in the central series of $A$. Thus $A_{(1)} = A$, and $A_{(n+1)} = [A_{(n)}, A]$.
Here we always consider closed subgroups of $A$. Hence $A_{(n+1)}$ is the closed subgroup of $A$ generated
by commutators $[x, y]$, $x \in A_{(n)}$, $y \in A$. For a $T$-group $A$ we define:

\begin{eqnarray*}
t_{1} & = & \dim_{{\fp}} H^{1}(\frac{Z(A)[p]}{T(A) \cap A_{(2)}}, \fp) \\
t_{i} & = & \dim_{{\fp}} H^{1}(\frac{Z(A)[p] \cap A_{(i)}}{T(A) \cap A_{(i+1)}}, \fp)  \  \  \  2 \le i \le p \\
\mu & = & \text{max}\{ i \colon 1 \le i\le p, A^{p} \subset A_{(i)} \}.
\end{eqnarray*}

These invariants are convenient for describing the $\fp[A/N]$-module $N$ associated with our $T$-group $A$;
see \cite[Section 1]{BLMS}. We have

\begin{prop}
For arbitrary cardinalities  $t_{i}$, $i =1, 2, \cdots, p$, and $\mu$ with $1 \le \mu \le p$, the following are
equivalent:
\begin{enumerate}
\item The $t_{i}$ and $\mu$ are invariants of some $T$-group.
\item  \begin{enumerate}
\item If $\mu < p$, then $t_{{\mu}} \ge 1$, and
\item If $\mu = p$ and $t_{i } = 0$ for all $2 \le i \le p$, then $t_{1} \ge 1$.
\end{enumerate}
\end{enumerate}
 Moreover, $T$-groups are uniquely determined up to isomorphism  by these invariants.
\end{prop}

 \begin{thm}
 For $p$ an odd prime, the following are equivalent.
 \begin{enumerate}
     \item $A$ is a $T$-group with invariants $t_{i}$ and $\mu$ satisfying
                \begin{enumerate}
                        \item $\mu \in \{1, 2\}$
                        \item $t_{2} = \mu -1$, and
                        \item $t_{i} = 0$ for $3 \le i \le p$.
                \end{enumerate}
     \item $A \cong T_{{E/F}}$ for some cyclic extension $E/F$ of degree $p$ such that
     $F$ contains a primitive $p$-th root of unity.
 \end{enumerate}
 Now suppose $p = 2$. Then each $T$-group is isomorphic to $T_{E/F}$ for some cyclic extension of
 degree $2$.
 \end{thm}

 It is interesting to note that these strong restrictions on possible $T$-groups occuring as $T_{{E/F}}$
 are consequences of the classical Hilbert 90 theorem which can now be viewed as the Bloch-Kato conjecture
 in degree $1$, as it  just involves basic Kummer theory which depends on Hilbert 90. For the realisation of
 given $T$-groups with invariants described in our theorem one uses constructions of cyclic extensions
 $E/F$ with prescribed groups $F^{*}/{F^{*}}^{p}$ and $N_{E/F}(E^{*})/{F^{*}}^{p}$ developed in \cite{MS05},
 which in turn uses results in \cite{Efrat-Maram} realising  certain semi-direct products and free
 pro-$p$-products as absolute Galois groups. This theorem provides restrictions on possible relations in $G_{F}(p)$; see
 \cite{BLMS}.

 \section{Higher Galois cohomology and the Bloch-Kato conjecture.}

In \cite{MS05} it was shown that the classical theorem Hilbert 90 is the key
for determining the $\fp[\Gal(E/F)]$-module structure of $E^{*}/{E^{*}}^{p}$
in the case of cyclic extensions of degree $p^{n}$. But
 $E^{*}/{E^{*}}^{p}$ is also $K_{1}(E)/pK_{1}(E)$. On the other hand Merkurjev-Suslin
 established in  \cite{Mer-Sus} an analogue   of the Hilbert 90 theorem for Milnor $K$-theory in degree $2$. Further it turned out that the analogue of Hilbert 90 for higher Milnor $K$-theory
 is essentially equivalent to the Bloch-Kato conjecture. This follows from the work of Merkurjev, Suslin, Rost and Voevodsky.
 Therefore it was  a natural idea to extend results on Galois module structure of $p$-power classes to Galois module structure
 of $K_{n}(E)/pK_{n}(E)$ for any positive integer $n$.  This was achieved in the case when $E/F$ is  a
 cyclic extension of prime
 degree $p$ and the primitive $p$-th root of unity is in $F$ in  \cite{LMS1}.
 (An extension of this work for the case of cyclic extension of degree $p^{n}$ is work in progress.)
 Some of  the main results are  explained in this section on the
 language of Galois cohomology as we freely use the Bloch-Kato conjecture.

Let $G$ be the Galois group of $E/F$. Write $E = F(\root{p}\of{a}), a \in
F^\times$ and denote $(a).(\xi_p) \in H^2(F,\fp)$ to be the cup product of
$(a),(\xi_p) \in H^1(F,\fp)$. For each $n \in \mathbb{N}$ set also:
  \begin{equation*}
    \Upsilon_1 \colon \dim_{\fp} \left( \ann_{H^{n-1}(F,\fp)}
    (a).(\xi_p)/\ann(a)\right)
  \end{equation*}
and
  \begin{equation*}
    \Upsilon_2 \colon \dim_{\fp} H^{n-1} (F,\fp)/\ann_{H^{n-1}
    (F,\fp)} (a).(\xi_p).
  \end{equation*}
Here $H^i(F,\fp) = H^i (G_F,\fp)$ is the $i$th Galois cohomology and $G_F$ is
the absolute Galois group of $F$. Here $\ann$ is an abbreviation for the
annihilator. Thus for example $\ann_{H^{n-1}(F,\fp)} (a).(\xi_p)$ is the
kernel of the cup product
  \begin{equation*}
    (a).(\xi_p).-\colon H^{n-1}(F,\fp) \to H^{n+1} (F,\fp).
  \end{equation*}
Set $U$ to be the absolute Galois group of $E$ and consider
$H^n(U,\fp) = H^n(E,\fp)$ as an $\fp[G]$-module. In \cite{LMS1}
we prove the following theorem. Our corestriction map here
denotes the map $H^n(E,\fp) \to H^n(F,\fp)$ and $\res$
means the restriction map $\res \colon H^n(F,\fp) \to
H^n(E,\fp)$.

\begin{thm} If $p>2$ and $n \in \mathbb{N}$ then
  \begin{equation*}
    H^n(E,\fp) \cong X_1 \oplus X_2 \oplus Y \oplus Z,
  \end{equation*}
where
  \begin{enumerate}
    \item $X_1$ is a trivial $\fp[G]$-module of dimension
    $\Upsilon_1$, and
      \begin{equation*}
        X_1 \cap \res H^n(F,\fp) = \{ 0 \}.
      \end{equation*}
    \item $X_2$ is a direct sum of $\Upsilon_2$ cyclic
    $\fp[G]$-modules of dimension $2$.
    \item $Y$ is a free $\fp[G]$-module of rank
      \begin{equation*}
        \dim_{\fp}  \text{Im}\left (\cor \colon \, H^n(E,\fp) \to H^n (F,\fp) \right)/
        (a).H^{n-1} (F,\fp) .
      \end{equation*}
    \item $Z$ is a trivial $\fp[G]$-module of dimension
      \begin{equation*}
        z = \dim_{\fp} H^n (F,\fp)/((\xi_p) H^{n-1} (F,\fp)
        + \cor\, H^n (E,\fp)) \mbox{ and } Z \subset \res
        (H^n(F,\fp)).
      \end{equation*}
    \end{enumerate}
\end{thm}

We see in particular that there is no cyclic summand of
dimension larger than $2$ but smaller than $p$.

A similar, but different theorem is valid in the case when
$p=2$. (See \cite[Theorem~2]{LMS1}.) Our decomposition of $G$-modules
$H^n(F,\fp)$ are not canonical but they allow a canonical
equivalent reformulation. Let $I$ be the augmentation ideal
of $\fp[G]$. Then from the analysis of the proof in Theorem~1,
one sees that our theorem is equivalent to the statements
below. (Here we abbreviate $H^n(E,\fp)$ as $H^n (E)$.)
  \begin{enumerate}
    \item For each $3 \le i \le p, I^{i-1} H^n (E) \cap
    H^n (E)^G = I^{p-1} H^n (E)$. \item $IH^m (E) \cap
    H^m (E)^G = \res (\xi_p).H^{m-1} (F) + \res \cor
    H^m (E)$. \item $0 \to \ann(a) \to H^{n-1} (E) \stackrel
    {(a).-}{\longrightarrow} H^m (F) \stackrel{\res}
    {\longrightarrow} H^m (E)^G \stackrel{\cor}{\longrightarrow}
    (a).\ann(a,\xi_p) \to 0$.
  \end{enumerate}
As we mentioned earlier this work uses in an essential way the key ingredients of the proof
of the Bloch-Kato conjecture, namely Hilbert~90 and the exact
sequence in Milnor $K$-theory. However, in the case of $n=2$,
we do not need the recent proof but we can use instead results from the
mid-1980s: the seminal results of Merkurjev and Suslin. (See
\cite{Mer-Sus}.)

In the case when characteristic of $F$ is $p$ and $G = \Gal(E/F) = \mathbb{Z}/p^{n}\mathbb{Z}$ is a 
cyclic group of prime power order all Galois modules $K_mE/p^s K_mE$, $s = 1,2, \cdots$ over
$\mathbb{Z}/p^{s}\mathbb{Z}[G]$ were classified; see \cite{BaLMS, MSS2}.

These results and ideas were applied to the development of an analog of
Schreier's formulas for the dimension of Galois cohomology under
$p$-extensions in  \cite{LLMS1}. In \cite{LLMS2} applications to the
characterization of Galois Demu\v{s}kin groups via Galois modules were
obtained. (Demu\v{s}kin groups are Poincar\'{e} groups of cohomological
dimension two, and Galois Demu\v{s}kin groups are Poincar\'{e} groups which
are also Galois groups of maximal $p$-extensions.)

Recall now that for a pro-$p$-group $G$ with finite
cohomology groups $H^i(G,\fp)$ for $0 \le i \le n$,
the $n$th partial Euler-Poincar\'e characteristic
$\chi_n(G)$ is defined as
  \begin{equation*}
    \chi_n(G) = \sum_{i=0}^{n} (-1)^i h_i (G).
  \end{equation*}
Suppose now that $G = G(p)$ for a field $F$ containing
a primitive $p$th root of unity, and suppose $G$ has
finite rank. Having determined $\chi_n(G)$, we have
obtained a strengthening of a theorem of Koch \cite[Theorem~5.5]{K}:

\noindent {\bf Theorem 3.} \cite[Corollary~2]{LMS2} \textit{Suppose that
$\xi_p \in F$, and let $n \in \mathbb{N}$. The following are equivalent:}
  \begin{enumerate}
    \item $\text{cd}(G) \le n$. \item $\chi_n(N) = p \chi_n(G)$ {\emph{for all
    open subgroups $N$ of $G$ of index}} $p$. \item $\chi_n(V) = p \chi_n
    (U)$ {\emph{for all open subgroups $U$ of $G$ and all open subgroups
    $V$ of $U$ of index}} $p$.
  \end{enumerate}

\section{Galois theoretic connections}
In this section we will explain the role played by certain Galois groups called W-groups in arithmetic.
To set the stage, we begin with our notation.

Let $F$ be a field of characteristic not equal to $2$. We shall introduce several subextensions
of $F_{sep}$.
\begin{itemize}
\item $F^{(2)}$ = compositum of all quadratic extensions of $F$.
\item  $F^{(3)}$ = compositum of all quadratic extensions of $F^{(2)}$, which are Galois over $F$.
%\item $F^{\{3\}}$ = compositum of all quadratic extensions of $F^{(2)}$.
\item $F_q$ = compositum of all Galois extensions $K/F$ such that $[K: F] = 2^n$, for some positive integer $n$.
\end{itemize}
All of these subextensions are Galois and they fit in a tower
\[ F \subset F^{(2)} \subset F^{(3)} \subset F_q \subset F_{sep} .\]
We denote their Galois groups (over $F$) as
\[G_F \lrar G_q \lrar G_F^{[3]} (=\G_F) \lrar G_F^{[2]} (= E) \lrar 1.\]

Observe that $G_F^{[2]}$ is just $\underset{i \in I}\prod C_2$, where $I$ is
the dimension of $F^*/{F^*}^2$ over $\ftwo$. $G_F$ is the absolute Galois
group of $F$. Although the quotients $G_q$ are much simpler than $G_{F}$ we
are far from understanding their structure in general. $F^{[3]}$ and its
Galois group over $F$ are considerably much simpler and yet they already
contain substantial arithmetic information of the absolute Galois group. The
groups $G_F^{[3]} (=\G_F)$ are called W-groups.

To illustrate this point, consider $WF$ the Witt ring of quadratic forms; see
\cite{Lam} for the definition. Then we have the following theorem.

\begin{thm}\cite{minac-spira-annals}
Let $F$ and $L$ be two fields of characteristic not $2$. Then $WF \cong WL$ (as rings) implies
that $\G_F \cong \G_L$ as pro-$2$-groups. Further if we assume additionally in the case when
each element of $F$ is a sum of two squares that $\sqrt{-1} \in F$ if and only if $\sqrt{-1} \in L$,
then $\G_F \cong \G_L$ implies $WF \cong WL$.
\end{thm}

Thus we see that $\G_F$ essentially controls the Witt ring $WF$ and in fact,
$\G_F$ can be viewed as a Galois theoretic analogue of $WF$. In particular,
$\G_F$ detects orderings of fields. (Recall that $P$ is an ordering of $F$ if
$P$ is an additively closed subgroup of index $2$ in  $F^*$.) More precisely,
we have:

\begin{thm} \cite{minac-spira-mathz}
There is a 1-1 correspondence between the orderings of a field $F$ and cosets
$\{ \sigma \Phi(\G_F) \, |\, \sigma \, \in \, \G_F \backslash \Phi(\G_F)  \text{and }\sigma^2 = 1\}$. Here $\Phi(\G_F)$
is the Frattini subgroup of $\G_F$, which is just the closed subgroup of $\G_F$ generated
by all squares in $\G_F$. The correspondence is as follows:
\[ \sigma \Phi(\G_F) \lrar P_{\sigma} = \{ f \in F^* \, | \, \sigma(\sqrt{f}) = \sqrt{f} \}.\]
\end{thm}

This theorem was generalised considerably for detecting additive properties of multiplicative
subgroups of $F^*$ in \cite{Louis-Minac-Smith}. In this paper (see  \cite[Section 8]{Louis-Minac-Smith})
it was shown that $\G_F$ can be used also for detecting valuations on $F$. The work on extending these ideas is in progress; see \cite{CEM}.

Also in \cite[Corollary 3.9]{adem-dikran-minac} it is shown that $\G_F \cong
\G_L$ if and only if $k_*(F) \cong k_*(L)$. Here $k_*(A)$ denotes the Milnor
$K$-theory (mod 2) of a field $A$. In particular, in  \cite[Theorem
3.14]{adem-dikran-minac} it is shown that if $R$ is the subring of $H^*(\G_F,
\ftwo)$ generated by one dimensional classes, then $R$ is isomorphic to the
Galois cohomology $ H^*(G_F, \ftwo)$ of $F$. Thus we see that $\G_F$ also
controls Galois cohomology and in fact $H^*(\G_F, \ftwo)$ contains some
further substantial information about $F$ which $H^*(G_F, \ftwo)$ does not
contain.  These results can be extended to the case when $p > 2$ and $F^{(2)}$ contains a
primitive $p$-th root of unity; see \cite{CEM, BCMS, BLMS}.
In summary, $\G_F$ is a very interesting object. On the one hand
$\G_F$ is much simpler than $G_F$ or $G_q$, yet it contains substantial
information about the arithmetic of $F$.   In fact, consider the case when  $p >2$ and $F$ contains a primitive
$p$-th root of unity. Then let $G = G_{F}$ be the absolute Galois group of $F$. The
descending $p$-central series of $G$ is defined inductively by  $G^{(1)} = G$, and
$G^{(i+1)} = (G^{(i)})^{p} [G^{(i)}, G]$, for $i \ge 1.$ Thus $G^{(i+1)}$ is the closed subgroup
of $G$ generated by all powers $h^{p}$ and all commutators $[h, g] = h^{-1}g^{-1}hg$, where
$h \in G^{(i)}$ and $g \in G$. Then the fixed fields $F^{(i)}$ of $G^{(i)}$ are precisely analogue
of fields
\[ F = F^{(1)}  \subset F^{(2)}  \subset F^{(3)}   \subset \cdots \subset F^{(i)} \subset \cdots \]
introduced above in the case $p =2$ and $i = 1, 2$ and $3$.

The special case of the main theorem in \cite{EM} then states:

\begin{thm} For $p > 2$ and for $G = G_{F}$ as above, $G^{(3)}$ is the intersection of
all normal subgroups $N$ of $G$ such that $G/N$ is isomorphic to one of $\{ 1\}$, $C_{p^{2}}$, and $M_{p^{3}}$ (the modular group of order $p^{3}$ which is the unique non-abelian group of
order $p^{3}$ and exponent $p^{2}$).
\end{thm}

The analogous result in the case $p =2$ was discovered by Villegas \cite{Vil}
in a different formalism. In \cite[Corollary 2.18]{minac-spira-annals} this
result was reformulated and reproved using the descending $2$-central sequence
of $G_{F}$. Namely, then $G^{(3)} = G_{F}^{(3)}$ is the intersection of all
open normal subgroups $N$ of $G$ such that $G/N$ is isomorphic to $\{ 1\},
C_{2}, C_{4}$, or to the dihedral group of order $8$.  The main ingredients in
the proofs of the above results is the Bloch-Kato conjecture in degree $2$
which was proved in \cite{Mer-Sus}.   The case $p =2$ was the first break-through
in the case of general fields made by Merkurjev who used some $K$-theoretic
calculations due to Suslin. For this particular case there is now an
elementary proof available due to Merkurjev; see  \cite{EKM}. If $p > 2$, in
the cohomology group $H^{2}(C_{p^{2}}, \mathbb{F}_{p})$ we have elements not
expressible as sums of products of elements in  $H^{1}(C_{p^{2}},
\mathbb{F}_{p})$. To handle these elements, in \cite{EM} there is a detailed
consideration of the Bockstein homomorphism
\[ B_{G} \colon H^{1}(G, \mathbb{F}_{p}) \rar H^{2}(G, \mathbb{F}_{p}).\]

In fact, in \cite{EM} not the full strength of Merkurjev-Suslin theorem was
used. The essential tool was the injectivity of the map
\[ k_{2} \rar H^{2}(F, \mathbb{F}_{p}).\]
In \cite{CEM}, the surjectivity of this map is used to obtain restrictions on presentation of
groups $G_{F}(p)$ via generators and relations.

Let $1 \rar R \rar S \rar G \rar 1,$
where $G = G_{F}(p)$ with $F$ as above, $S$ a free pro-$p$-group with minimal number
of generators (see \cite[Chapter 4]{K}), and $R$ is the subgroup of $S$ of relations in $G$.
Then we have:

\begin{thm} (\cite{CEM} for any $p$, \cite{minac-spira-annals} for $p =2$.)
Let $S \supset S^{(2)} \supset S^{(3)} \supset \cdots$ be the $p$-descending series of $S$. Then we have
\[R^{p} [R, S] = R \bigcap S^{(3)}.\]
\end{thm}

Observe that for any minimal presentation of any pro-$p$-group $G$ as above, we have
$R^{p} [R, S] \subset R \bigcap S^{(3)}$, as $R \subset S^{(2)}$. The equality in the case
when $G = G_{F}(p)$ is the consequence of the surjectivity of the norm residue map
\[ k_{2}(F) \rar H^{2}(G_{F}(p), \mathbb{F}_{p}),\]
which follows from the Merkurjev-Suslin's theorem.

From the above theorem, one can deduce that if $G$ is any pro-$p$-group such
that $R \subset S^{(3)}$, and $G = G_{F}(p)$, then $G$ is a free
pro-$p$-group.

\medskip

\noindent
\textbf{Example} Let $G$ be a pro-$p$-group on $n$ generators $a_{1}, a_{2}, \cdots, a_{n}$
for $n \ge 2$ subject to relations $[[a_{i}, a_{j}], a_{r}]  = 1$ for all $1 \le i < j \le n$ and $1 \le r \le n$.  Then $G$ is not $G_{F}(p)$ for any field $F$ containing a primitive $p$th root of unity. \cite{CEM} contains further restrictions on possible groups $G_{F}(p)$ by exploring further
properties of its quotients $G_{F}^{[3]} = G_{F}(p)/G_{F}(p)^{(3)}$ and close connections
between properties of $G_{F}^{[3]}$ and the existence of non-trivial valuations on $F$.

We now outline a joint project with Benson and Swallow in which our goal is to  obtain a refinement of the Bloch-Kato conjecture.  Associated to the field $F$, we have a
natural tower of subfields $F^{(n)}$ of the separable closure $F_{\sep}$ defined as follows: $F^{(1)} = F$,
$F^{(2)}$ is the  compositum of all cyclic extensions of degree $p$ over $F$, and for $n \ge 3$, $F^{(n)}$ is
the compositum of all cyclic extensions of degree $p$ over $F^{(n+1)}$ which are Galois over $F$. We call this tower
the \emph{filtration tower} associated to $F$. We define Galois groups (in agreement with the notation we used above)
\begin{align*}
G_F^{[n]} & := \Gal(F^{(n)}/F),  \ \ \text{and}\\
G_F^{(n)} & := \Gal(F_{\sep}/F^{(n)}),
\end{align*}
which fit in a sequence $1 \rar G_F^{(n)} \rar G_F\rar G_F^{[n]}\rar 1$, where
$G_F$ is the absolute Galois group $\Gal(F_{\sep}/F)$.  In \cite{CEM} (see
also \cite{BCMS}),  we have shown that the decomposable part of
$H^*(G_{F}^{[3]}, \mathbb{F}_p)$ is isomorphic to  $H^*(F , \mathbb{F}_p)$
under the inflation map. The important question therefore is to determine how
an indecomposable class in $H^*(G_{F}^{[3]}, \mathbb{F}_p)$ decomposes under
the various inflation maps along the filtration tower. By the Bloch-Kato
conjecture, we know that it decomposes completely into one-dimensional classes
when it goes all the way up to the separable closure. But what happens in
between? A precise knowledge of this gives a refinement of the Bloch-Kato
conjecture. We have shown (using the Bloch-Kato conjecture in degree 2!) that
every indecomposable class in $H^2(G_F^{[n]}, \mathbb{F}_p)$ decomposes into
one-dimensional classes when it goes to the next level $H^2(G_F^{[n+1]},
\mathbb{F}_p)$ under the inflation map. Thus we have obtained a second
cohomology refinement of the Bloch-Kato conjecture. The goal of our joint
project with Benson and Swallow is to understand this refinement of the
Bloch-Kato conjecture for higher cohomology. This is work in progress.

\section*{Acknowledgements}
We are very grateful to our collaborators and friends especially to Dave
Benson, Ido Efrat,  John Labute, Andy Schultz, and John Swallow for working
with us  on these fascinating topics and sharing with us their unique insight.
We would like also to thank the referee, whose valuable suggestions helped us
to improve our exposition.

%\bibliographystyle{alpha}
%\bibliography{lit}

\begin{thebibliography}{AWKDM}

\bibitem[AKM99]{adem-dikran-minac}
Alejandro Adem, Dikran~B. Karagueuzian, and J{\'a}n Min{\'a}{\v{c}}.
\newblock On the cohomology of {G}alois groups determined by {W}itt rings.
\newblock {\em Adv. Math.}, 148(1):105--160, 1999.

\bibitem[AGKM]{AAWKDM}
Alejandro Adem, Wenfeng Gao, Karagueuzian Dikran B, and  J{\'a}n Min{\'a}{\v{c}}.
\newblock Field theory and the cohomology of some Galois groups.
\newblock J. Algebra  235  (2001),  no. 2, 608--635.

\bibitem[BCMS]{BCMS}
D.~J. Benson, Sunil~K. Chebolu, J\'{a}n Min\'{a}\v{c}, and John Swallow.
\newblock Bloch-Kato pro-p groups and a refinement of the Bloch-Kato
  conjecture.
\newblock 2007.
\newblock preprint.

\bibitem[BLMS]{BLMS}
D. J. Benson, Nicole Lemire, J{\'a}n Min{\'a}{\v{c}}, and John Swallow.
\newblock Detecting pro-{$p$}-groups that are not absolute {G}alois groups.
\newblock {\em J. Reine Angew. Math.}, 613:175--191, 2007.

\bibitem[BhLMS]{BaLMS}
Ganesh Bhandari, Nicole Lemire, J{\'a}n Min{\'a}{\v{c}}, and John Swallow.
\newblock Galois module structure of Milnor $K$-theory in characteristic $p$.
\newblock {\em New York J. Math.}, 14:215-224, 2008.

\bibitem[Bor65]{Borevic}
Z.~I. Borevi{\v{c}}.
\newblock The multiplicative group of cyclic {$p$}-extensions of a local field.
\newblock {\em Trudy Mat. Inst. Steklov}, 80:16--29, 1965.

\bibitem[CEM]{CEM}
 Sunil~K. Chebolu, Ido Efrat, and J\'{a}n Min\'{a}\v{c}.
\newblock The structure of Galois groups of maximal pro-$p$ quotients of absolute Galois groups.
\newblock 2008.
\newblock Preprint.


\bibitem[DuLa]{DuLa} D. Dummit and J. Labute. On a new characterization
of Demu\v{s}kin groups. \emph{Invent. Math.} {\bf 73} (1983), no. 3,
413--418.

\bibitem[EKM]{EKM} R. Elman, N. Karpenko, and A. Merkurjev.
\newblock The algebraic and geometric theory
of quadratic forms,
\newblock  AMS, Colloquium publications, Vol 56, 2008.

\bibitem[EfHa]{Efrat-Maram}
I. Efrat and D. Haran.
\newblock On Galois groups over Pythagorean and semi-real closed fields.
\newblock {\em Israel Journal of Math.}, 85:57-87, 1996.

\bibitem[EM]{EM}
 Ido Efrat, and J\'{a}n Min\'{a}\v{c}.
\newblock   On the descending central sequence of absolute Galois groups
\newblock Preprint 2008.
\newblock  arXiv:0809.2166

\bibitem[Fad60]{Faddeev}
D.~K. Faddeev.
\newblock On the structure of the reduced multiplicative group of a cyclic
  extension of a local field.
\newblock {\em Izv. Akad. Nauk SSSR Ser. Mat.}, 24:145--152, 1960.

\bibitem[HW]{HW} Ch. Haesemeyer and C. W. Weibel. Norm varieties
and the chain lemma (after Markus Rost). Proc. Abel Symposium (to appear), 31
pages, 2008 preprint.

\bibitem[K]{K} H. Koch. \emph{Galois theory of $p$-extensions}. Berlin:
Springer-Verlag, 2002.

\bibitem[LLMS06]{LLMS2} J. Labute, N. Lemire, J. Min\'a\v{c}, and J.
Swallow. Demu\v{s}kin groups, Galois modules, and the elementary
type conjecture. \emph{J. Algebra}, 304:1130--1146, 2006.

\bibitem[LLMS07]{LLMS1} J. Labute, N. Lemire, J. Min\'a\v{c}, and J.
Swallow. Cohomological dimension and Schreier's formula in Galois
cohomology. \emph{Can. Math. Bull.},  50: 588-593, 2007.

\bibitem[Lam05]{Lam}
T.~Y. Lam.
\newblock {\em Introduction to quadratic forms over fields}, volume~67 of {\em
  Graduate Studies in Mathematics}.
\newblock American Mathematical Society, Providence, RI, 2005.

\bibitem[LMS05]{LMS2} N.~Lemire, J.~Min\'a\v{c}, and J.~Swallow.
When is Galois cohomology free or trivial? \emph{New York J.~Math.} {\bf 11} (2005), 291--302.

\bibitem[LMS07]{LMS1} N. Lemire, J. Min\'a\v{c}, and J. Swallow.
Galois module structure of Galois cohomology and partial
Euler-Poincar\'e characteristics.
\newblock {\em J. Reine Angew. Math}, 613: 147 -- 173, (2007).

\bibitem[MMS04]{Louis-Minac-Smith}
Louis Mah{\'e}, J{\'a}n Min{\'a}{\v{c}}, and Tara~L. Smith.
\newblock Additive structure of multiplicative subgroups of fields and {G}alois
  theory.
\newblock {\em Doc. Math.}, 9:301--355 (electronic), 2004.

\bibitem[Mer81]{Mer81}
Merkurjev, A. S.
\newblock On the norm residue symbol of degree $2$. (Russian)
\newblock Dokl. Akad. Nauk SSSR  261  (1981), no. 3, 542--547.

\bibitem[MS82]{Mer-Sus}
A.~S. Merkurjev and A.~A. Suslin.
\newblock {$K$}-cohomology of {S}everi-{B}rauer varieties and the norm residue
  homomorphism.
\newblock {\em Izv. Akad. Nauk SSSR Ser. Mat.}, 46(5):1011--1046, 1135--1136,
  1982. (Translation from Russian: \emph{Math. USSR-Izv.} {\bf 21}
(1983), 307--340.)


\bibitem[Mil70]{Mil-1970}
John Milnor.
\newblock Algebraic {$K$}-theory and quadratic forms.
\newblock {\em Invent. Math.}, 9:318--344, 1969/1970.

\bibitem[MS90]{minac-spira-mathz}
J{\'a}n Min{\'a}{\v{c}} and Michel Spira.
\newblock Formally real fields, {P}ythagorean fields, {$C$}-fields and
  {$W$}-groups.
\newblock {\em Math. Z.}, 205(4):519--530, 1990.

\bibitem[MS96]{minac-spira-annals}
J{\'a}n Min\'{a}\v{c} and Michel Spira.
\newblock Witt rings and {G}alois groups.
\newblock {\em Ann. of Math. (2)}, 144(1):35--60, 1996.
\bibitem[MVW]{MVW} C.~Mazza, V.~Voevodsky, and C.~W.~Weibel.
\emph{Lecture notes on motivic cohomology}. Clay Mathematics Monographs 2.
Providence, RI: American Mathematical Society; Cambridge, MA: Clay
Mathematics Institute, 2006.

\bibitem[MS03]{minac-swallow}
J{\'a}n Min{\'a}{\v{c}} and John Swallow.
\newblock Galois module structure of {$p$}th-power classes of extensions of
  degree {$p$}.
\newblock {\em Israel J. Math.}, 138:29--42, 2003.

\bibitem[MS05]{MS05}
J{\'a}n Min{\'a}{\v{c}} and John Swallow.
\newblock Galois modules  appearing as {$p$}th-power classes of units of extensions of
  degree {$p$}.
\newblock {\em Math. Zeit.}, 250, no.4, 2005.

\bibitem[MSS06]{MSS}
J{\'a}n Min{\'a}{\v{c}}, Andrew Schultz, and John Swallow.
\newblock Galois module structure of {$p$}th-power classes of cyclic extensions
  of degree {$p\sp n$}.
\newblock {\em Proc. London Math. Soc. (3)}, 92(2):307--341, 2006.

\bibitem[MSS2]{MSS2} Jan  Min\'a\v{c}, A. Schultz, and J. Swallow.
Galois module structure of Milnor $K$-theory mod $p^s$ in
characteristic $p$.
\newblock {\em New York J. Math}, 14: 225 -- 233, (2008).

\bibitem[R1]{Ro} M.~Rost. Chain lemma for symbols.
Available at www.math.uni-bielefeld.de/$\sim$rost/chain-lemma.html.

\bibitem[R2]{Ro2} M.~Rost.  On the basic
correspondence of a splitting variety. Available at
www.math.uni-bielefeld.de/$\sim$rost/chain-lemma.html.

\bibitem[SJ]{Su} A.~A.~Suslin and S.~Joukhovitski. Norm
varieties. {\it J.~Pure Appl.~Algebra} {\bf 206} (2006),
235--276.

\bibitem[Ser02]{serre-GaloisCohomology}
Jean-Pierre Serre.
\newblock {\em Galois cohomology}.
\newblock Springer Monographs in Mathematics. Springer-Verlag, Berlin, English
  edition, 2002.
\newblock Translated from the French by Patrick Ion and revised by the author.


\bibitem[Vi]{Vil} F. R. Villegas.
Relations between quadratic forms and certain Galois extensions, a manuscript,
Ohio State University, 1988,
http://www.math.utexas.edu/users/villegas/osu.pdf.

\bibitem[V1]{Vo1} V.~Voevodsky. Motivic cohomology with
$\mathbb{Z}/2$-coefficients. \emph{Publ. Inst. Hautes \'Etudes Sci.} {\bf 98}
(2003), 59--104.

\bibitem[V2]{Vo2} V.~Voevodsky. On motivic cohomology with
$\mathbb{Z}/l$-coefficients. $K$-theory preprint archive 639.
Available at www.math.uiuc.edu/K-theory/0639/.

\bibitem[V3]{Vo3} V.~Voevodsky. Motivic Eilenberg-MacLane
spaces. $K$-theory preprint archive 864.
Available at www.math.uiuc.edu/K-theory/0864/.

\bibitem[W1]{W1} C.~W.~Weibel. The norm residue isomorphism
theorem, 18 pp, 2007 preprint. Available at
www.math.rutgers.edu/$\sim$weibel/papers.html.

\bibitem[W2]{W2} C.~W.~Weibel. Axioms for the norm residue
isomorphism. pp. 427-435 in K-theory and Noncommutative Geometry, European
Math. Soc. Pub. House, 2008

\bibitem[W3]{W3} C.~W.~Weibel. The proof of the Bloch-Kato
Conjecture. ICTP Lecture Notes Series 23 (2008), 1-28


\end{thebibliography}

\end{document}